\documentclass[10pt,conference]{IEEEtran}
\usepackage{amsmath,epsfig}
\usepackage[defs]{ams}
\usepackage[usenames]{color}
\usepackage{umoline}
\usepackage{setspace}
\usepackage{algorithm}
\usepackage{algorithmic}
\usepackage{multirow}

\newcommand{\INPUT}{\item[\textbf{Input:}]}
\newcommand{\OUTPUT}{\item[\textbf{Output:}]}

\newcommand{\norm}[1]{\left\Vert#1\right\Vert}

\newcommand{\supp}[1]{\mathrm{supp}(#1)}

\newcommand{\rank}{\mathrm{rank}}
\newcommand{\atoms}[1]{\mathrm{atoms}(#1)}

\newtheorem{theorem}{Theorem}[section]

\newtheorem{proposition}[theorem]{Proposition}

\newtheorem{definition}[theorem]{Definition}
\newtheorem{remark}[theorem]{Remark}

\newcommand{\qed}{\nobreak \ifvmode \relax \else
      \ifdim\lastskip<1.5em \hskip-\lastskip
      \hskip1.5em plus0em minus0.5em \fi \nobreak
      \vrule height0.75em width0.5em depth0.25em\fi}

\begin{document}

\title{Efficient and Guaranteed Rank Minimization by Atomic Decomposition \vspace{-4mm}}

\author{
\authorblockN{Kiryung Lee and Yoram Bresler}
\authorblockA{Coordinated Science Laboratory and Department of ECE\\
University of Illinois at Urbana-Champaign\\
1308 W. Main St., Urbana, IL 61801\\
Email: \{klee81,ybresler\}@illinois.edu}
\vspace{-8mm}
}

\maketitle

\begin{abstract}
Recht, Fazel, and Parrilo provided an analogy between rank minimization and $\ell_0$-norm minimization.
Subject to the rank-restricted isometry property,
nuclear norm minimization is a guaranteed algorithm for rank minimization.
The resulting semidefinite formulation is a convex problem but in practice the algorithms for it do not scale well to large instances.
Instead, we explore missing terms in the analogy and propose a new algorithm which is computationally efficient and also has a performance guarantee.
The algorithm is based on the atomic decomposition of the matrix variable and extends the idea in the CoSaMP algorithm for $\ell_0$-norm minimization.
Combined with the recent fast low rank approximation of matrices based on randomization,
the proposed algorithm can efficiently handle large scale rank minimization problems.
\end{abstract}

\section{Introduction}
\label{sec:intro}

Recent studies in compressed sensing have shown that
a sparsity prior in the representation of the unknowns can guarantee unique and stable solutions to underdetermined linear systems.
The idea has been generalized to linear rank minimization by Recht, Fazel, and Parrilo \cite{recht2007gmr}.
Rank minimization has important applications such as
matrix completion, linear system identification, Euclidean embedding, and image compression.

The rank minimization problem is formally written as:
\begin{equation*}
\text{P1:} \qquad
\begin{array}{llll}
\displaystyle \min_{X \in \mathbb{C}^{m \times n}} & \rank(X) \\
\mathrm{subject~to} & \mathcal{A}X = b,
\end{array}
\end{equation*}
for a given linear operator $\mathcal{A}: \mathbb{C}^{m \times n} \rightarrow \mathbb{C}^{p}$ and $b \in \mathbb{C}^p$.

Fazel, Hindi, and Boyd \cite{fazel2001rmh} proposed a convex relaxation of the rank minimization problem.
They minimized the nuclear norm $\norm{X}_*$, which is the sum of all singular values of matrix $X$,
and is the convex envelop of the non-convex function $\rank(X)$.
Rank minimization is related to $\ell_0$-norm minimization, which has been the focus of compressed sensing.
Recht, Fazel, and Parrilo \cite{recht2007gmr} provided an analogy
between the two problems and their respective solutions by convex relaxation.

In the analogy, $\ell_1$-norm minimization for the $\ell_0$-norm minimization problem is replaced by nuclear norm minimization for $\text{P1}$.
Both are efficient algorithms, with guaranteed performance under certain conditions, to solve NP-hard problems:
$\ell_0$-norm minimization and rank minimization, respectively.
The respective conditions are given by the restricted isometry property and its generalization.
However, whereas $\ell_1$-norm minimization corresponds to a linear program,
nuclear norm minimization is formulated as a convex semidefinite program (SDP).
Although there exist polynomial time algorithms to solve SDP, in practice they do not scale well to large problems.

Recently, several authors proposed methods for solving large scale SDP derived from rank minimization.
These include interior point methods for SDP, projected subgradient methods, and low-rank parametrization \cite{recht2007gmr}
and a customized interior point method \cite{liu2008ipm}  These methods can solve larger rank minimization problems,
which the general purpose SDP solvers cannot.
However, the dimension of the problem is still restricted and some of these methods do not guarantee convergence to the global minimum.
Other methods solve nuclear norm minimization in
a penalized form using singular value thresholding (SVT) \cite{cai2008svt} or fixed point and Bregman iterations \cite{ma2008fpb}.
It has been shown that the sequence of solutions converges to the solution to nuclear norm minimization as the penalty parameter increases.
However, an analysis of the convergence rate is missing and hence the quality of the solution obtained by these methods is not guaranteed.
Furthermore, the efficiency of these methods is restricted to the case of an affine (i.e., linear equality) constraint.
Meka \textit{et. al.} \cite{meka2008rmv} used multiplicative updates and online convex programming to provide an approximate solution to rank minimization.
However, their result depends on the (unverified) existence of an oracle that provides the solution to the rank minimization problem with a single linear constraint in constant time.

For $\ell_0$-norm minimization, besides $\ell_1$-norm minimization,
there are recent algorithms, which are more efficient and also have performance guarantees.
These include Compressive Sampling Matching Pursuit (CoSaMP) \cite{tropp2008cis} and Subspace Pursuit (SP) \cite{dai2008spc}.
To date, no such algorithms have been available for rank minimization.

In this paper, we propose an iterative algorithm to solve the rank minimization problem,
which is a generalization
\footnote{
There is another generalization of CoSaMP, namely model-based CoSaMP \cite{baraniuk2008mbc}.
However, this generalization addresses a completely different and unrelated problem: sparse vector approximation subject to a special (e.g., tree) structure.
Furthermore, the extensions of CoSaMP to model-based CoSaMP and to ADMiRA are independent:
neither one follows from the other, and neither one is a special case of the other.
}
of the CoSaMP algorithm for $\ell_0$-norm minimization to the rank minimization problem.
We call this algorithm ``Atomic Decomposition for Minimum Rank Approximation,'' abbreviated as ADMiRA.
In CoSaMP, the $\ell_0$-norm minimization problem with equality constraints is recast into an $s$-term vector approximation problem.
Similarly, in ADMiRA we recast $\text{P1}$ into the rank-$r$ matrix approximation problem
\begin{equation*}
\text{P2:} \qquad
\begin{array}{llll}
\displaystyle \min_{X \in \mathbb{C}^{m \times n}} & \norm{\mathcal{A}X - b}_2 \\
\mathrm{subject~to} & \rank(X) \leq r.
\end{array}
\end{equation*}
The performance guarantee of ADMiRA states that
the approximate solution to $\text{P2}$ obtained by ADMiRA coincides with the true solution $X$ to $\text{P1}$
for any $r \geq \rank(X)$ that satisfies the ADMiRA assumptions.

\section{Vector vs Matrix}
\label{sec:vector_vs_matrix}

\subsection{Preliminaries}

Throughout this paper, we use two vector spaces:
the space of column vectors $\mathbb{C}^p$ and the space of matrices $\mathbb{C}^{m \times n}$.
For $\mathbb{C}^p$, the inner product is defined by
$\langle x , y \rangle_{\mathbb{C}^p} = y^H x$ for $x, y \in \mathbb{C}^p$ where $y^H$ denotes the Hermitian transpose of $y$,
and the induced Hilbert-Schmidt norm is the Euclidean or $\ell_2$-norm given by
$\norm{x}_2^2 = \langle x, x \rangle_{\mathbb{C}^p}$ for $x \in \mathbb{C}^p$.

For $\mathbb{C}^{m \times n}$, the inner product is defined by
$\langle X, Y \rangle_{\mathbb{C}^{m \times n}} = \mathrm{tr}(Y^H X)$ for $X, Y \in \mathbb{C}^{m \times n}$,
and the induced norm is the Frobenius norm given by
$\norm{X}_F^2 = \langle X, X \rangle_{\mathbb{C}^{m \times n}}$ for $X \in \mathbb{C}^{m \times n}$.

\subsection{Atomic Decomposition}

Let $S$ denote the set of all nonzero rank-one matrices in $\mathbb{C}^{m \times n}$.
We can refine $S$ so that any two distinct elements are not collinear.
The resulting subset $\mathbb{O}$ is referred to as the \textit{set of atoms}
\footnote{The ``atom'' in this paper is different from Mallat and Zhang's ``atom'' \cite{mallat1993mpt},
which is an element in the dictionary, a finite set of vectors.
In both cases, however, an atom denotes an irreducible quantity.
}
of $\mathbb{C}^{m \times n}$.
Then the \textit{set of atomic spaces} $\mathbb{A}$ of $\mathbb{C}^{m \times n}$ is defined by
\begin{equation}
\mathbb{A} \triangleq \{\mathrm{span}(\psi) :~ \psi \in \mathbb{O}\}.
\end{equation}
Each subspace $V \in \mathbb{A}$ is one-dimensional and hence is irreducible in the sense that
$V = V_1 + V_2$ for some $V_1,V_2 \in \mathbb{A}$ implies $V_1 = V_2 = V$.
Since $\mathbb{O}$ is a uncountably infinite set in a finite dimensional space $\mathbb{C}^{m \times n}$,
the elements in $\mathbb{O}$ are not linearly independent.
Regardless of the choice of $\mathbb{O}$, $\mathbb{A}$ is uniquely determined.
Without loss of generality, we fix $\mathbb{O}$ such that all elements have the unit Frobenius norm.

Given a matrix $X \in \mathbb{C}^{m \times n}$,
its representation $X = \sum_j \alpha_j \psi_j$ as a linear combination of atoms
is referred to as an \textit{atomic decomposition} of $X$.
Since $\mathbb{O}$ spans $\mathbb{C}^{m \times n}$, an atomic decomposition of $X$ exists for all $X \in \mathbb{C}^{m \times n}$.
A subset $\Psi = \{ \psi \in \mathbb{O} :~ \langle \psi_j, \psi_k \rangle_{\mathbb{C}^{m \times n}} = \delta_{jk} \}$
of unit-norm and pairwise orthogonal atoms in $\mathbb{O}$ will be called an \textit{orthonormal set of atoms}.

\begin{definition}
Let $\mathbb{O}$ be a set of atoms of $\mathbb{C}^{m \times n}$.
Given $X \in \mathbb{C}^{m \times n}$, we define $\atoms{X}$
\begin{equation}
\atoms{X} \triangleq \displaystyle \arg\min_\Psi \left\{ |\Psi| :~ \Psi \subset \mathbb{O}, \quad X \in \mathrm{span}(\Psi) \right\}.
\label{eq:def_atom}
\end{equation}
\end{definition}
Note that $\atoms{X}$ is not unique.

An orthonormal set $\atoms{X}$ is given by the singular value decomposition of $X$.
Let $X = \sum_{k=1}^{\rank(X)} \sigma_k u_k v_k^H$ denote the singular value decomposition of $X$ with singular values in decreasing order.
For each $k$, there exists $\rho_k \in \mathbb{C}$ such that $|\rho_k| = 1$ and $\rho_k u_k v_k^H \in \mathbb{O}$.
Then an orthonormal set $\atoms{X}$ is given by
\begin{equation*}
\atoms{X} = \{ \rho_k u_k v_k^H \}_{k=1}^{\rank(X)}.
\end{equation*}

\begin{remark}
$\atoms{X}$ and $\rank(X) = |\atoms{X}|$ of a matrix $X \in \mathbb{C}^{m \times n}$ are the counterparts of
$\supp{x}$ and $\norm{x}_0 = |\supp{x}|$ for a vector $x \in \mathbb{C}^p$, respectively.
\end{remark}

\subsection{Generalized Correlation Maximization}

Recht, Fazel, and Parrilo \cite{recht2007gmr} showed an analogy between rank minimization $\text{P1}$ and $\ell_0$-norm minimization.
We consider instead the rank-$r$ matrix approximation problem $\text{P2}$ and its analogue -- the $s$-term vector approximation problem
\begin{equation*}
\text{P3:} \qquad
\begin{array}{llll}
\displaystyle \min_{x \in \mathbb{C}^{n}} & \norm{Ax - b}_2 \\
\mathrm{subject~to} & \norm{x}_0 \leq s.
\end{array}
\end{equation*}
In Problem $\text{P3}$, variable $x$ lives in the union of $s$ dimensional subspaces of $\mathbb{C}^{n}$,
each spanned by $s$ elements in the finite set $\mathbb{E} = \{ e_1, \ldots, e_n \}$, the standard basis of $\mathbb{C}^n$.
Thus the union contains all $s$-sparse vectors in $\mathbb{C}^n$.
Importantly, finitely many ($n \choose s$, to be precise) subspaces participate in the union.
Therefore, it is not surprising that $\text{P3}$ can be solved exactly by exhaustive enumeration,
and finite selection algorithms such as CoSaMP are applicable.

In the rank-$r$ matrix approximation problem $\text{P2}$,
the matrix variable $X$ lives in the union of subspaces of $\mathbb{C}^{m \times n}$,
each of which is spanned by $r$ atoms in the set $\mathbb{O}$.
Indeed, if $X \in \mathbb{C}^{m \times n}$ is spanned by $r$ atoms in $\mathbb{O}$, then
$\rank(X) \leq r$ by the subadditivity of the rank.
Conversely, if $\rank(X) = r$,
then $X$ is a linear combination of rank-one matrices and hence there exist $r$ atoms that span $X$.
Note that uncountably infinitely many subspaces participate in the union.
Therefore, some selection rules in the greedy algorithms for $\ell_0$-norm minimization and $s$-term vector approximation do not generalize in a straightforward way.
None the less, using our formulation of the rank-$r$ matrix approximation problem in terms of an atomic decomposition,
we extend the analogy between the vector and matrix cases,
and propose a way to generalize these selection rules to the rank-$r$ matrix approximation problems.

First, consider the correlation maximization in greedy algorithms for the vector case.
Matching Pursuit (MP) \cite{mallat1993mpt} and Orthogonal Matching Pursuit (OMP) \cite{pati1993omp} choose the index $k \in \{1,\ldots,n\}$
that maximizes the correlation $\left|a_k^H (b - A \hat{x})\right|$ between the $k$-th column $a_k$ of $A$ and the residual in each iteration,
where $\hat{x}$ is the solution of the previous iteration.
Given a set $\Psi$, let $P_\Psi$ denote the projection operator onto the subspace spanned by $\Psi$ in the corresponding embedding space.
When $\Psi = \{ \psi \}$ is a singleton set, $P_\psi$ will denote $P_\Psi$.
For example, $P_{e_k}$ denotes the projection operator onto the subspace in $\mathbb{C}^n$ spanned by $e_k$.
From
\begin{eqnarray*}
\left| a_k^H (b - A \hat{x}) \right|
= \left| \langle A^H(b - A \hat{x}), e_k \rangle_{\mathbb{C}^n} \right|
= \norm{P_{e_k} A^H (b - A \hat{x})}_2,
\end{eqnarray*}
it follows that
maximizing the correlation implies
maximizing the norm of the projection of the image under $A^H$ of the residual $b - A \hat{x}$ onto the selected one dimensional subspace.

The following selection rule generalizes the correlation maximization to the matrix case.
We maximize the norm of the projection over all one-dimensional subspaces spanned by an atom in $\mathbb{O}$:
\begin{equation}
\max_{\psi \in \mathbb{O}} \left| \langle b - \mathcal{A} \widehat{X}, \mathcal{A} \psi \rangle_{\mathbb{C}^{m \times n}} \right|
= \max_{\psi \in \mathbb{O}} \norm{P_\psi \mathcal{A}^*(b - \mathcal{A} \widehat{X})}_F,
\label{eq:sel_one_largest_cor}
\end{equation}
where $\mathcal{A}^* : \mathbb{C}^p \rightarrow \mathbb{C}^{m \times n}$ denotes the adjoint operator of $\mathcal{A}$.
By the Eckart-Young Theorem, the basis of the best subspace is obtained from the singular value decomposition of $M = \mathcal{A}^*(b - \mathcal{A} \widehat{X})$,
as $\psi = u_1 v_1^H$, where $u_1$ and $v_1$ are the principal left and right singular vectors.
\begin{remark}
Applying the selection rule (\ref{eq:sel_one_largest_cor}) to update $\widehat{X}$ recursively leads to greedy algorithms generalizing MP and OMP to rank minimization.
\end{remark}

Next, consider the rule in recent algorithms such as CoSaMP, and SP.
The selection rule chooses the subset $J$ of $\{1,\ldots,n\}$ with $|J| = s$ defined by
\begin{equation}
\left|a_k^H (b - A \hat{x})\right| \geq \left|a_j^H (b - A \hat{x})\right|, \quad \forall k \in J, \forall j \not\in J.
\label{eq:sel_largest_cor}
\end{equation}
This is equivalent to maximizing
\begin{eqnarray*}
\sum_{k \in J} \left|a_k^H (b - A \hat{x})\right|^2
&=& \norm{P_{\{e_k\}_{k \in J}} A^H (b - A \hat{x})}_2^2.
\end{eqnarray*}
In other words, selection rule (\ref{eq:sel_largest_cor}) finds the best subspace spanned by $s$ elements in $\mathbb{E}$
that maximizes the norm of the projection of $M = A^H (b - A \hat{x})$ onto that $s$-dimensional subspace.

The following selection rule generalizes the selection rule (\ref{eq:sel_largest_cor}) to the matrix case.
We maximize the norm of the projection over all subspaces spanned by a subset with at most $r$ atoms in $\mathbb{O}$:
\begin{equation*}
\max_{\Psi \subset \mathbb{O}} \left\{ \norm{P_\Psi \mathcal{A}^*(b - \mathcal{A} \widehat{X})}_F :~ |\Psi| \leq r \right\}
\end{equation*}
A basis $\Psi$ of the best subspace is again obtained from the singular value decomposition of $M = \mathcal{A}^*(b - \mathcal{A} \widehat{X})$,
as $\Psi = \{\rho_k u_k v_k^H\}_{k=1}^r$, where $u_k$ and $v_k$, $k = 1,\ldots,r$ are the $r$ principal left and right singular vectors, respectively
and for each $k$, $\rho_k \in \mathbb{C}$ satisfies $|\rho_k| = 1$
\footnote{
Once the best subspace is determined, it is not required to compute the constants $\rho_k$'s.
}
.
Note that $\Psi$ is an orthonormal set although this is not enforced as an explicit constraint in the maximization.

\section{Algorithm}
\label{sec:algorithm}

\begin{algorithm}
\caption{ADMiRA}
\begin{algorithmic}[1]
\INPUT $\mathcal{A}: \mathbb{C}^{m \times n} \rightarrow \mathbb{C}^p$, $b \in \mathbb{C}^p$, and target rank $r \in \mathbb{N}$
\OUTPUT rank-$r$ solution $\widehat{X}$ to $\text{P2}$
\STATE $\widehat{X} \leftarrow 0$
\STATE $\widehat{\Psi} \leftarrow \emptyset$
\WHILE{stop criterion is false}
\STATE $\Psi' \leftarrow \displaystyle \arg\max_{\Psi \subset \mathbb{O}} \left\{ \norm{P_\Psi \mathcal{A}^* (b - \mathcal{A}\widehat{X})}_F :~ |\Psi| \leq 2r \right\}$
\label{step:cor_max_A}
\STATE $\widetilde{\Psi} \leftarrow \Psi' \cup \widehat{\Psi}$
\label{step:merge}
\STATE $\widetilde{X} \leftarrow \displaystyle \arg\min_X \left\{ \norm{b - \mathcal{A}X}_2 :~ X \in \mathrm{span}(\widetilde{\Psi}) \right\}$
\label{step:solve_LS}
\STATE $\widehat{\Psi} \leftarrow \displaystyle \arg\max_{\Psi \subset \mathbb{O}} \left\{ \norm{P_\Psi \widetilde{X}}_F :~ |\Psi| \leq r \right\}$
\label{step:cor_max_C}
\STATE $\widehat{X} \leftarrow P_{\widehat{\Psi}} \widetilde{X}$
\label{step:proj}
\ENDWHILE
\RETURN $\widehat{X}$
\end{algorithmic}
\label{alg:ADMiRA}
\end{algorithm}

Algorithm~\ref{alg:ADMiRA} describes ADMiRA.
Steps~\ref{step:cor_max_A} and \ref{step:cor_max_C} involve finding a best rank-$2r$ or rank-$r$ approximation to given matrix (e.g., by truncating the SVD),
while Step~\ref{step:solve_LS} involves the solution of a linear least-squares problem -- all standard numerical linear algebra problems.
Step~\ref{step:merge} merges two given sets of atoms in $\mathbb{O}$ by taking their union.

Most steps of ADMiRA are similar to those of CoSaMP except Step~\ref{step:cor_max_A} and Step~\ref{step:cor_max_C}.
The feasible sets of the maximization problems in Step~\ref{step:cor_max_A} and Step~\ref{step:cor_max_C} of ADMiRA are infinite
while those in the analogous steps of CoSaMP are finite.
In CoSaMP, a greedy selection is employed to solve the combinatorial problem and provides the exact solution owing to the orthogonality of the feasible set.
The maximization problem over the infinite set in ADMiRA may look even more difficult than the combinatorial problem in CoSaMP.
However, singular value decomposition can solve the maximization problem over the infinite set efficiently.

\section{Main Results: Performance Guarantee}

\subsection{Rank-Restricted Isometry Property (R-RIP)}

Recht \textit{et al} \cite{recht2007gmr} generalized the sparsity-restricted isometry property (RIP) defined for sparse vectors to low rank matrices.
\footnote{They also demonstrated ``nearly isometric families'' satisfying this R-RIP (with overwhelming probability).
These include random linear operators generated from i.i.d. Gaussian or i.i.d. symmetric Bernoulli distributions.}
In order to draw the analogy with known results in $\ell_0$-norm minimization, we slightly modify their definition by squaring the norm in the inequality.
Given a linear operator $\mathcal{A}: \mathbb{C}^{m \times n} \rightarrow \mathbb{C}^p$,
the rank-restricted isometry constant $\delta_r(\mathcal{A})$ is defined as the minimum constant that satisfies
\begin{equation}
(1 - \delta_r(\mathcal{A})) \norm{X}_F^2 \leq \norm{\gamma \mathcal{A} X}_2^2 \leq (1 + \delta_r(\mathcal{A})) \norm{X}_F^2,
\label{eq:rip}
\end{equation}
for all $X \in \mathbb{C}^{m \times n}$ with $\rank(X) \leq r$ for some constant $\gamma > 0$.
Throughout this paper, we assume that the linear operator $\mathcal{A}$ is scaled appropriately so that $\gamma = 1$ in (\ref{eq:rip})
\footnote{
If $\gamma \neq 1$, then only the constant for the noise gain will be scaled accordingly.
}
.

\subsection{Performance Guarantee}

Subject to the R-RIP, the Atomic Decomposition for Minimum Rank Approximation Algorithm (ADMiRA) has a performance guarantee analogous to that of CoSaMP.

The followings are the assumptions in ADMiRA:
\begin{description}
\item[A1:] The target rank is fixed as $r$.
\item[A2:] The linear operator $\mathcal{A}$ satisfies $\delta_{7r}(\mathcal{A}) \leq 0.043$.
\item[A3:] The measurement is obtained by
\begin{equation}
b = \mathcal{A}X + \nu,
\label{eq:noisy_meas_mdl}
\end{equation}
where $\nu$ is the discrepancy between the measurement and the linear model $\mathcal{A}X$.
\end{description}

Assumption A2 plays a key role in deriving the performance guarantee of ADMiRA.
This enforces the rank-restricted isometry property of the linear operator $\mathcal{A}$.
Although the verification of the satisfiability of A2 is as difficult as or more difficult than the recovery problem itself,
nearly isometric families that satisfy the condition in A2 have been demonstrated \cite{recht2007gmr}.

The performance guarantees are specified in terms of a measure of inherent approximation error, termed \emph{unrecoverable energy} defined by
\begin{equation*}
\epsilon = \norm{X - X_r}_F + \frac{1}{\sqrt{r}} \norm{X - X_r}_* + \norm{\nu}_2,
\end{equation*}
where $X_r$ denotes the best rank-$r$ approximation of $X$.
The first two terms in $\epsilon$ define a metric of the minimum distance
between the ``true'' matrix $X$ and a rank-$r$ matrix.
This is analogous to the notion of a measure of compressibility of a vector in sparse vector approximation.
No solution of $\text{P2}$ can come any closer to $X$.
The third term is the norm of the measurement noise, which must also limit the accuracy of the approximation provided by a solution to $\text{P2}$.

\begin{theorem}
Let $\widehat{X}_k$ denote the estimate of $X$ in the $k$-th iteration of ADMiRA.
For each $k \geq 0$, $\widehat{X}_k$ satisfies the following recursion:
\begin{equation*}
\|X - \widehat{X}_{k+1}\|_F \leq 0.5 \|X - \widehat{X}_k\|_F + 10 \epsilon,
\end{equation*}
where $\epsilon$ is the unrecoverable energy.
From the above relation, it follows that
\begin{equation*}
\|X - \widehat{X}_k\|_F \leq 2^{-k} \norm{X}_F + 20 \epsilon, \quad \forall k \geq 0.
\end{equation*}
\label{thm:pg_gen}
\end{theorem}
\vspace{-3mm}
Theorem~\ref{thm:pg_gen} shows the geometric convergence of ADMiRA.
In fact, convergence in a finite number of steps can be achieved as stated by the following theorem.

\begin{theorem}
After at most $6(r+1)$ iterations, ADMiRA provides a rank-$r$ approximation $\widehat{X}$ of $X$, which satisfies
\begin{equation*}
\|X - \widehat{X}\|_F \leq 20 \epsilon,
\end{equation*}
where $\epsilon$ is the unrecoverable energy.
\label{thm:iter_cnt_gen}
\end{theorem}
Depending on the spectral properties of the matrix $X$, even faster convergence is possible.

\subsection{Relationship between $\text{P1}$, $\text{P2}$, and ADMiRA}

The approximation $\widehat{X}$ given by ADMiRA is a solution to $\text{P2}$.
When there is no noise in the measurement, i.e., $b = \mathcal{A}X$, where $X$ is the solution to $\text{P1}$,
Theorem~\ref{thm:pg_gen} states that if the ADMiRA assumptions are satisfied with $r \geq \rank(X)$, then $\widehat{X} = X$.
An appropriate value can be assigned to $r$ by an incremental search over $r$.

For the noisy measurement case, the linear constraint in $\text{P1}$ is replaced by a quadratic constraint and
the rank minimization problem is written as:
\begin{equation*}
\text{P1$'$:} \qquad
\begin{array}{llll}
\displaystyle \min_{X \in \mathbb{C}^{m \times n}} & \rank(X) \\
\mathrm{subject~to} & \norm{\mathcal{A}X - b}_2 \leq \eta.
\end{array}
\end{equation*}
Let $X'$ denote a minimizer to $\text{P1$'$}$.
In this case, the approximation $\widehat{X}$ produced by ADMiRA is not necessarily equivalent to $X'$,
but by Theorem~\ref{thm:pg_gen} the distance between the two is bounded by $\|X' - \widehat{X}\|_F \leq 20 \eta$
for all $r \geq \rank(X')$ that satisfies the ADMiRA assumptions.

\section{Properties of the Rank-Restricted Isometry}

We introduce a number of properties of the rank-restricted isometry.
These properties serve as key tools for proving the performance guarantees for ADMiRA in this paper.
These properties further extend the analogy between the sparse vector and the low-rank matrix approximation problems
($\text{P3}$ and $\text{P2}$, respectively), and are therefore also of interest in their own right.
An operator satisfying the R-RIP satisfies, as a consequence, a number of other properties
when composed with other linear operators defined by the atomic decomposition.
Most properties are inherited from the vector case.
However, the generalization of the \emph{restricted orthogonality property} to the matrix case is not straightforward
and shows some nontrivial differences.
The following Proposition is an extension of Lemma~2.1 in \cite{candes2008rip} for the vector case to the matrix case.

\begin{proposition}
Suppose that linear operator $\mathcal{A}: \mathbb{C}^{m \times n} \rightarrow \mathbb{C}^p$ has the rank-restricted isometry constant $\delta_r(\mathcal{A})$.
Let $X, Y \in \mathbb{C}^{m \times n}$ such that $\langle X, Y \rangle_{\mathbb{C}^{m \times n}} = 0$ and $\rank(X) + \rank(Y) \leq r$.
Then
\begin{equation}
\left| \langle \mathcal{A} X , \mathcal{A} Y \rangle_{\mathbb{C}^p} \right| \leq \sqrt{2} \delta_r(\mathcal{A}) \norm{X}_F \norm{Y}_F.
\end{equation}
\label{prop:rop_mat}
\end{proposition}
\vspace{-3mm}

\begin{remark}
For the real matrix case, Proposition~\ref{prop:rop_mat} can be improved by dropping the constant $\sqrt{2}$.
This improvement is achieved by replacing the parallelogram identity in the proof to the version for the real scalar field case.
\end{remark}

\begin{remark}
For the vector case, the representation of a vector $x \in \mathbb{C}^n$ in terms of the standard basis $\{e_j\}_{j=1}^n$ of $\mathbb{C}^n$ determines $\norm{x}_0$.
Let $J_1, J_2 \subset \{1,\ldots,n\}$ be arbitrary.
Then the projection operators $P_{\{e_j\}_{j \in J_1}}$ and $P_{\{e_j\}_{j \in J_2}}$ commute.
Furthermore, $P_{\{e_j\}_{j \in J_1}}^\perp x$ is $s$-sparse (or sparser) if $x$ is $s$-sparse.
These properties follow from the orthogonality of the standard basis.
Proposition 3.2 in \cite{tropp2008cis}, corresponding in the vector case to our Proposition~\ref{prop:rop_mat} requires these two properties.
However, these properties do not hold for the matrix case.
For $\Psi_1, \Psi_2 \subset \mathbb{O}$, the projection operators $P_{\Psi_1}$ and $P_{\Psi_2}$ do not commute in general
and $\rank(P_\Psi X)$ can be greater than $r$ even though $\rank(X) \leq r$.
Proposition~\ref{prop:rop_mat} is a stronger version of the corresponding proposition for the vector case
in the sense that it requires a weaker condition (orthogonality between two low-rank matrices),
which can be satisfied without these properties.
\end{remark}

\section{Implementation and Scalability}
\label{sec:large}

Most of the computation in ADMiRA lies in the truncated singular value decomposition.
The fact that ADMiRA keeps the matrix variables in their atomic decomposition is advantageous for this procedure.
Only a few dominant singular triplets are necessary, which can be computed by the Lanczos method in $O((m+n)rL)$ time,
where $L$ is the number of the iterations that depends on the singular value distribution..
An alternative approach is to use a randomized algorithm \cite{harpeled2006lrm}
that computes the low-rank approximation of a given matrix in atomic decomposed form in $O((m+n)r^3\log r)$ time.
In this case, ADMiRA has complexity of $O((m+n)r^3\log r)$ per iteration, or $O((m+n)r^4\log r)$ to achieve the guarantee in Theorem~\ref{thm:iter_cnt_gen},
and scales well to large problems.
Another consideration is the computation of the proxy matrix.
This involves applying $\mathcal{A}$ and $\mathcal{A}^*$, the complexity of which is $O(rpmn)$.
If $\mathcal{A}$ consists of sparse matrices, then the complexity can be as small as $O(rp(m+n))$.
In particular, in the matrix completion problem,
$\mathcal{A}X$ is sampling the entries of matrix $X$ and hence there is no multiplication in this procedure.

\section{Numerical Experiment}
\label{sec:num_exp}

\begin{figure}[htb]
\begin{center}
\begin{minipage}[htb]{0.45\linewidth}
\centerline{\includegraphics[height=40mm]{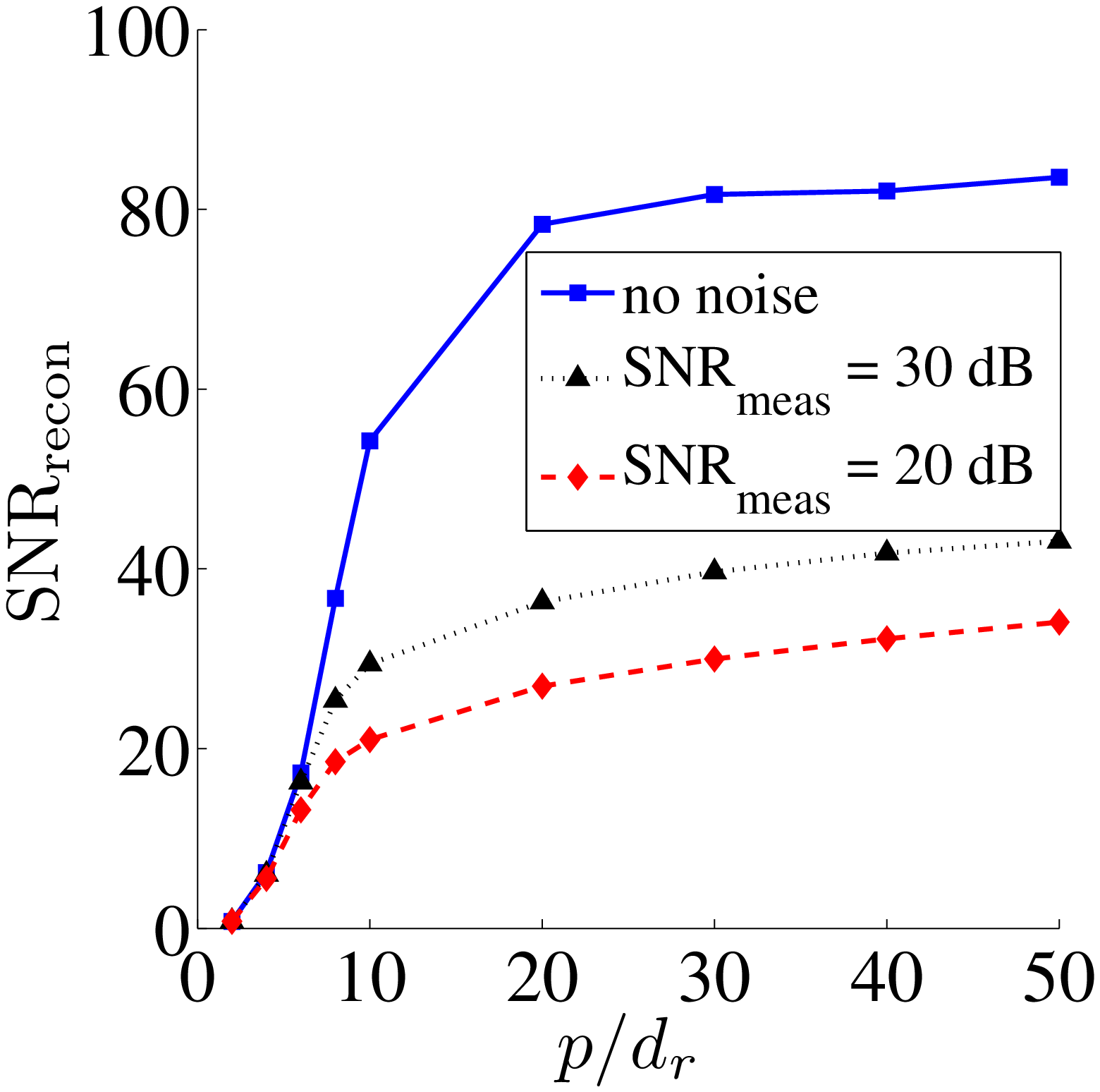}}
\end{minipage}
\begin{minipage}[htb]{0.45\linewidth}
\centerline{\includegraphics[height=40mm]{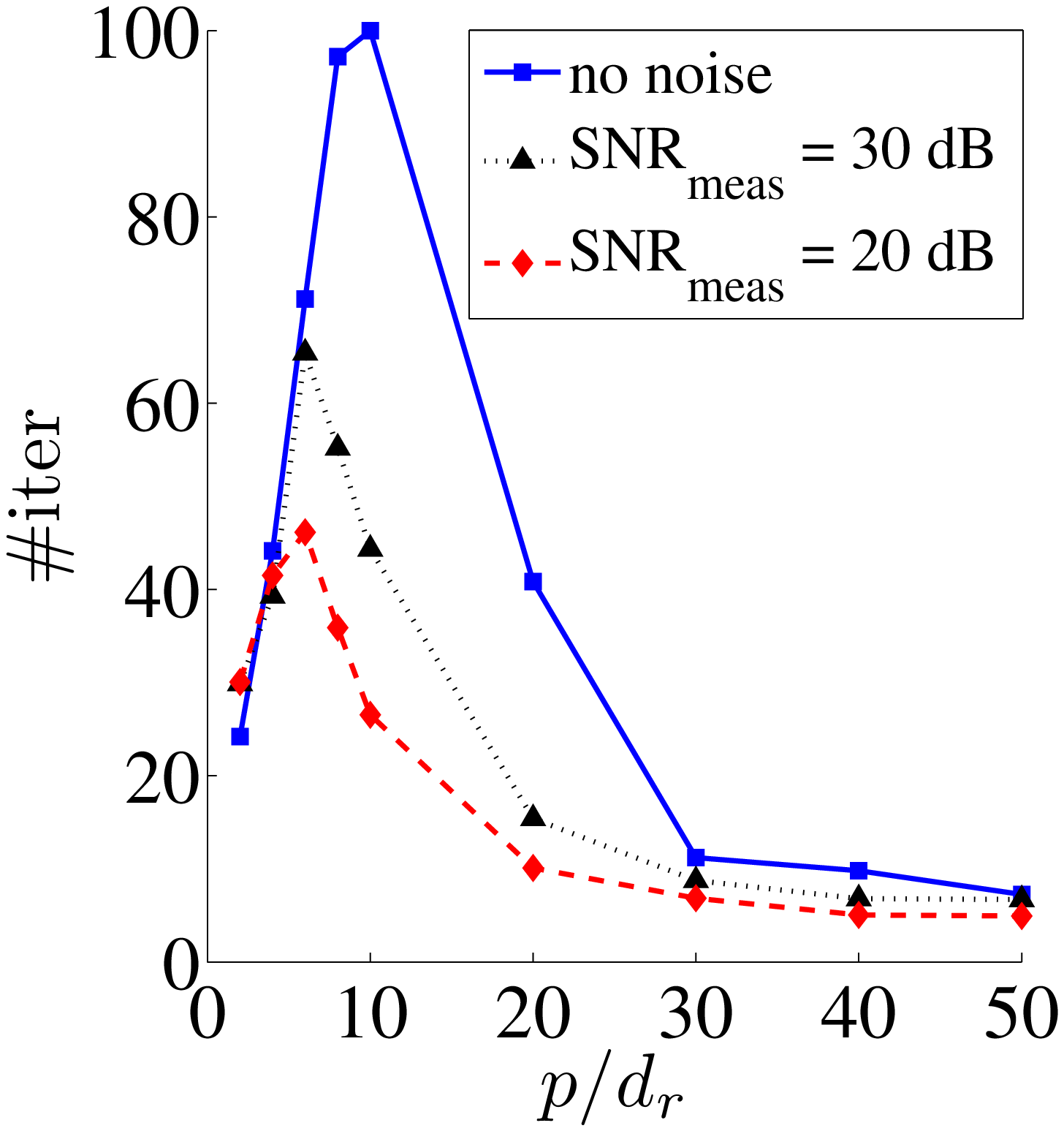}}
\end{minipage}
\end{center}
\vspace{-2mm}
\caption{Matrix completion by ADMiRA: $n = 500, r = 2$.}
\vspace{-2mm}
\label{fig:var_p}
\end{figure}

We study reconstructions by ADMiRA with a generic matrix completion example.
Our preliminary Matlab implementation uses ARPACK \cite{lehoucq:aug} to compute partial SVDs in Steps~\ref{step:cor_max_A} and \ref{step:cor_max_C} of ADMiRA.
The test matrix $X \in \mathbb{R}^{n \times n}$ is generated as the product $X = Y_L Y_R^H$
where $Y_L,Y_R \in \mathbb{R}^{n \times r}$ has entries following an i.i.d. Gaussian distribution.
The measurement $b$ is $p$ randomly chosen entries of $X$, which may be contaminated with an additive white Gaussian noise.
The reconstruction error and measurement noise level are measured in terms of
$\mathrm{SNR}_\mathrm{recon} \triangleq 20 \log_{10}(\norm{X}_F / \|X - \widehat{X}\|_F)$ and
$\mathrm{SNR}_\mathrm{meas} \triangleq 20 \log_{10}(\norm{b}_2 / \norm{\nu}_2)$, respectively.
Computational efficiency is measured by the number of iterations.
The results in Fig.~\ref{fig:var_p} and Table~\ref{table:compare_admira_svt} have been averaged over 20 trials.

Fig.~\ref{fig:var_p} shows that both $\mathrm{SNR}_\mathrm{recon}$ and the number of iterations improve as $p/d_r$ increases.
Here $d_r$ is the number of degrees of freedom defined by $d_r = r(n+m-r)$ and denotes the essential number of unknowns.
Fig.~\ref{fig:var_p} suggests that we need $p / d_r \geq 20$ for $n = 500$.

Table~\ref{table:compare_admira_svt} shows that ADMiRA provides slightly better performance with less computation than SVT \cite{cai2008svt}.
Roughly, the computational complexity of a single iteration of ADMiRA can be compared to three times of that of SVT.

\begin{figure}[htb]
\begin{center}
\begin{minipage}[htb]{0.45\linewidth}
\centerline{\includegraphics[height=40mm]{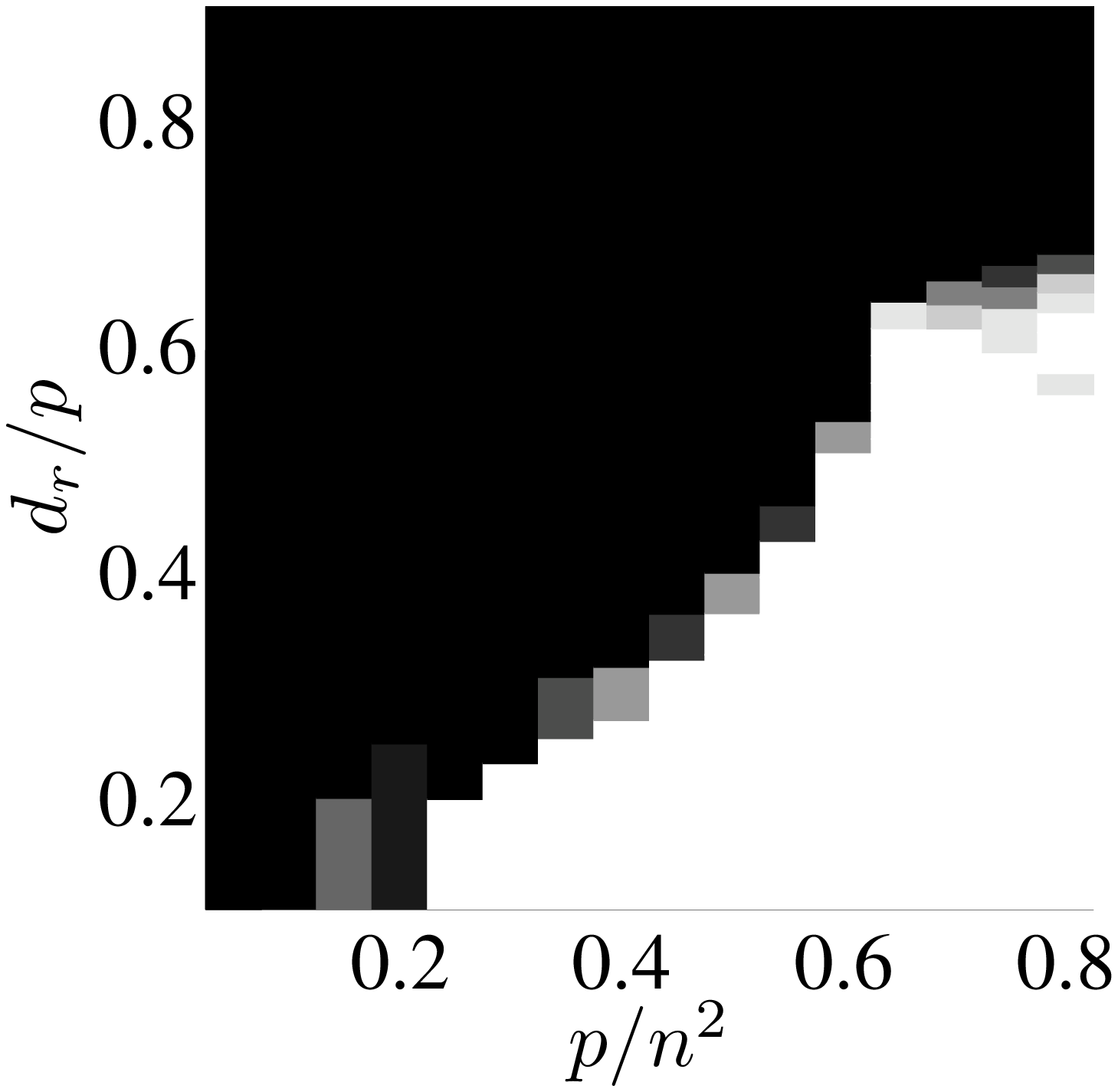}}
\vspace{-1mm}
\centering{\footnotesize ADMiRA}
\end{minipage}
\begin{minipage}[htb]{0.45\linewidth}
\centerline{\includegraphics[height=40mm]{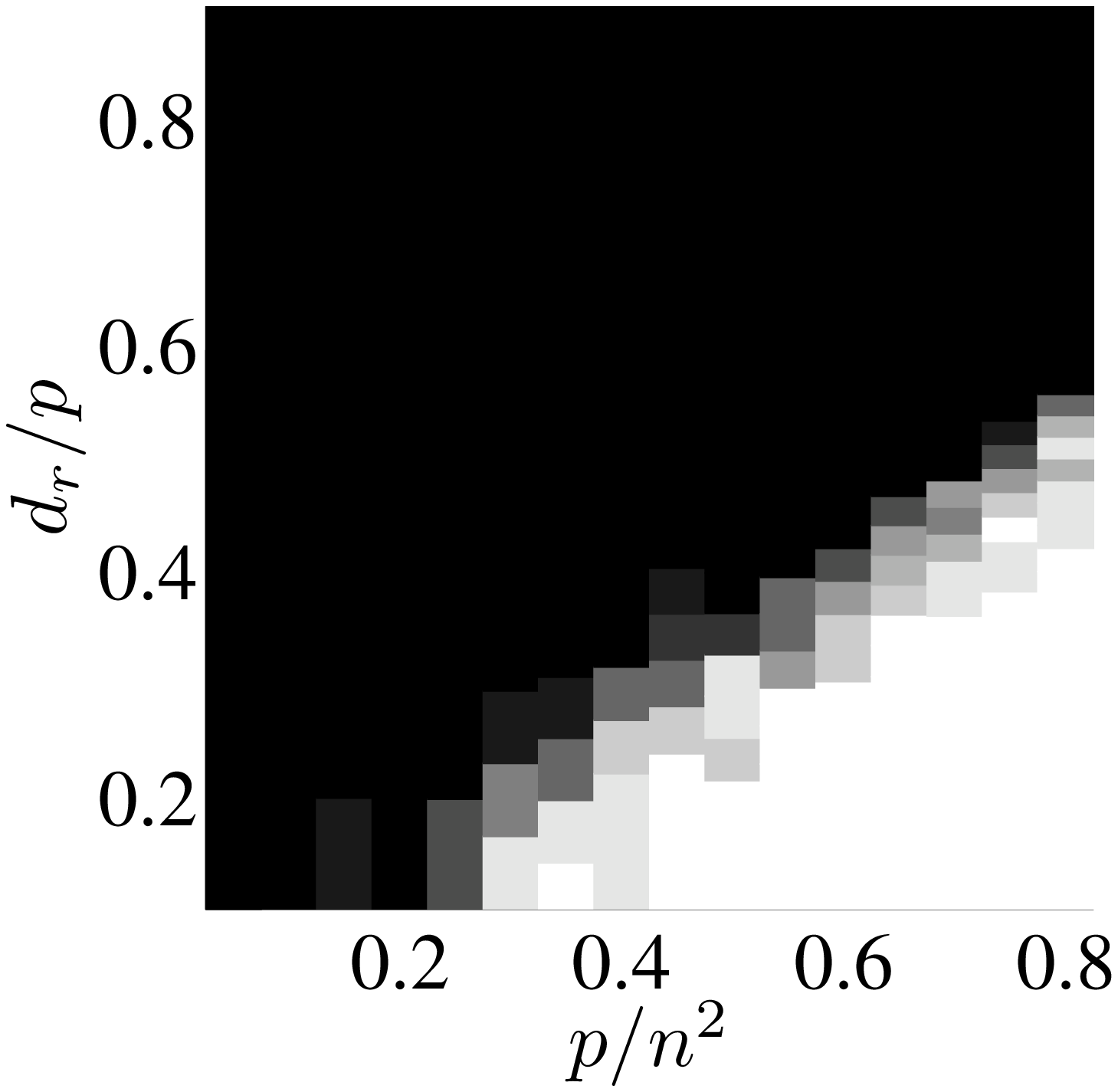}}
\vspace{-1mm}
\centering{\footnotesize SVT}
\end{minipage}
\end{center}
\vspace{-2mm}
\caption{Phase transition of matrix completion: $n = m = 100$.}
\label{fig:ptd}
\end{figure}

Fig.~\ref{fig:ptd} compares the phase transitions of ADMiRA and SVT.
We count the number of successful matrix completion ($\mathrm{SNR}_\mathrm{recon} \geq 70 \mathrm{dB}$) out of 10 trials for each triplet $(n,p,r)$.
The brighter color implies more success.
ADMiRA performed better than SVT for this example.

\begin{table}
\begin{center}
\footnotesize
\begin{minipage}[b]{0.9\linewidth}
\centerline{\begin{tabular}{c|*{5}{c|}c}
  \hline
  \multirow{2}{*}{$r$} & \multirow{2}{*}{$p/n^2$} & \multirow{2}{*}{$p/d_r$} & \multicolumn{2}{|c|}{$\mathrm{SNR}_\mathrm{recon}$ (dB)} & \multicolumn{2}{|c}{$\#\text{iter}$} \\\cline{4-7}
  & & & ADMiRA & SVT & ADMiRA & SVT \\\hline\hline
  2 & 0.20 & 50.05 & 82 & 79 & 11 & 54 \\\hline
  5 & 0.20 & 20.05 & 81 & 78 & 15 & 64 \\\hline
  10 & 0.20 & 10.05 & 79 & 77 & 19 & 80 \\\hline
\end{tabular}}
\end{minipage}
\end{center}
\vspace{-2mm}
\caption{Comparison of ADMiRA and SVT: no noise, $n = m = 1000$.}
\vspace{-6mm}
\label{table:compare_admira_svt}
\end{table}

We note though, that the performance guarantee in the previous sections is not directly applicable to the experiments in this section
since the linear operator in the matrix completion does not satisfy the R-RIP.
It seems that a performance guarantee without using the R-RIP might be possible.

\section{Conclusion}
\label{sec:conclusion}

We propose a new algorithm, ADMiRA,
which extends both the efficiency and the performance guarantee of the CoSaMP algorithm for $\ell_0$-norm minimization to matrix rank minimization.
The proposed generalized correlation maximization can be also applied to MP, OMP, and SP
to similarly extend the known algorithms and theory from the $s$-term vector approximation problem to the rank-$r$ matrix approximation.
ADMiRA can handle large scale rank minimization problems efficiently by using recent linear time algorithms for low rank approximation.
More detailed arguments and missing proofs are available in \cite{leebre2009admirafull}.
\vspace{-2mm}



\end{document}